\documentclass[9pt, twoside]{amsart}
\usepackage{amssymb}
\usepackage{float}
\usepackage{amsfonts}
\usepackage[centertags]{amsmath}
\usepackage{amsthm}
\usepackage{newlfont}
\usepackage{array}
\usepackage{graphicx, color}

\date{}

\begin{document}

\centerline{\Large{\bf New type Pythagorean fuzzy soft set}}

\centerline{}

\centerline{\Large{\bf and decision-making application}}

\centerline{}

\centerline{\bf {Murat Kiri\c{s}ci}}

\centerline{}

\centerline{Department of Mathematical Education, Hasan Ali Y\"{u}cel Education Faculty,}

\centerline{ Istanbul University-Cerrahpa\c{s}a, Vefa, 34470, Fatih, Istanbul, Turkey}

\centerline{}

\newtheorem{Theorem}{\quad Theorem}[section]

\newtheorem{Definition}[Theorem]{\quad Definition}

\newtheorem{Corollary}[Theorem]{\quad Corollary}

\newtheorem{Proposition}[Theorem]{\quad Proposition}

\newtheorem{Lemma}[Theorem]{\quad Lemma}

\newtheorem{Example}[Theorem]{\quad Example}

\newtheorem*{remark}{Remark}

\centerline{}

{\bf Abstract:} We define the Pythagorean fuzzy parameterized soft set and investigate some properties of the new set.
Further, we propose to the solution of decision-making application for the Pythagorean
fuzzy parameterized soft set and other related concepts.\\



\section{Introduction}

Many terms that we use randomly in everyday life usually have a fuzzy structure.
Verbal or numerical expressions we use, while describing something, explaining an event,
commanding and in many other cases include fuzziness.  People use terms that do not express
certainty when explaining an event and deciding on a situation.  According to the age of
the person, old, middle, young, very old and very young are called. Depending on the slope
and ramp condition of the road, the car's gas or brake pedal is press slightly slower or
slightly faster.  All these are examples of how the human brain behaves in uncertain and
vagueness  situations, and how it evaluates, identifies, and commands events.\\

After the fuzzy set(FS) theory, which uses fuzzy logic rules was developed by Lotfi A. Zadeh and
published in its original 1965 paper \cite{Zadeh}, the examination of uncertainty systems has gained a new dimension.
FSs are characterized by membership functions. In fact, these membership functions are nothing more
than fuzzy numbers. A set defined in this way can be described by a membership function that appoints membership values
to all of the elements from 0 to 1. Members that are not included in the set are appointed membership values
of 0, and those who are included in the set are appointed membership values of 1. The elements that are not
included in the set are appointed values between 0 and 1 according to uncertainty situation.\\

With the emergence of FS theory, science and technology have made great improvement.

It is seen to have many applications related to the FS Theory in both theoretical and practical studies from health sciences to computer science, from physical sciences to arts, and from engineering and humanities to life sciences.
For theoretical study examples, refer to the \cite{chang}, \cite{kirisci3}, \cite{matloka}, \cite{zararsiz0}, \cite{zararsiz1}, \cite{zararsiz2}.\\

The problems we face in our lives are often not clear and precise.
So we use various decision-making mechanisms to solve our problems.
With these mechanisms we use, we try to make the most right decision by reducing uncertainties.
Therefore, improved mathematical tools for uncertainty and imprecision are needed.
Soft Set Theory(SST) has been used fairly broadly to deal with such imprecision.\\

The concept of SS has been initiated by Molodtsov. Molodtsov proposed the concept of the SS, an entirely new approach to modeling uncertainty.
The SST has a generous application potential.
Some of these applications have been shown by Molodtsov in his pioneering work.
This theory was implemented in many areas of uncertainty such as mathematical analysis, algebraic structures, optimization theory, information systems, decision-making problems.
Maji and et al.\cite{MajiBis}
investigated the SS theory for decision-making problems.
It is investigated by the SS Theory to decision-making problems. \cite{MajiBis1}.
In \cite{MajiBis2}, the authors defined operations "AND", "OR", union and intersection of two SSs.
The same authors established a hybrid model known as fuzzy soft set(FSS), which is a unification of SS and FS \cite{MajiBis}.
Actually, FSS is an extension of crisp SS.\\

Recently, after the generalization of SSs as FSSs, different extensions have been included in the literature such as
intuitionistic FSSs \cite{MajiBis2}, fuzzy parametrized SSs \cite{CagEn}, intuitionistic fuzzy parametrized SSs \cite{DeCa}.
Some developments in SSs, FSs, and applications in decision-making can be found in \cite{Alietal}, \cite{Fengetal},
\cite{kirisci0}, \cite{kirisci1}, \cite{kirisci2}, \cite{Yager2}.\\

Yager \cite{Yager0} offered a new FS called Pythagorean fuzzy set(PFS).
PFS has fascinated the care of great deal researchers  in a little while time.
Yager and Abbasov \cite{YagerAbb} improved the concept of Pythagorean membership grades. The formulation of the negation for IFSs and PFSs is examined by Yager \cite{Yager1}.
In \cite{YagerKit},PF subsets and its relationship with IF subsets were debated and some set operations on PF subsets were defined.\\

In \cite{PengYang}, the properties such as boundedness, idempotency, and monotonicity related to the Pythagorean fuzzy aggregation operators are investigated. Further, for to solve uncertainty
multiple attribute group decision-making problem Pythagorean fuzzy superiority and inferiority ranking method was developed in \cite{PengYang}.\\

Peng et al.\cite{Pengetal}, defined the PFSS and investigated its properties.
Guleria and Bajaj \cite{Guleria} proposed PF soft matrix and its diverse feasible types.
Additionally, the PF soft matrices have been well-considered for recommending a new algorithm for decision-making by using choice matrix and weighted choice matrix.\\

We introduce Pythagorean fuzzy parameterized soft set and investigate some properties, operations.
Further, we present to the solution of decision-making problem with Pythagorean fuzzy parameterized soft set and other related concepts.\\

In section 2, we give some basic definitions and properties. In section 3, we define the Pythagorean fuzzy
parameterized soft set and study basic properties. In section 4, we give basic operations according to Pythagorean fuzzy parameterized soft set
In section 5, we introduce the results related to decision-making Problem.

\section{Preliminaries}

The FS has emerged as a generalization of the classical set concept.
If we choose a non-empty set $X$, then a function $m_{A}(x):X\rightarrow [0,1]$
is called FS on $X$ and represented by
\begin{eqnarray*}
A=\left\{(x_{i},m_{A}(x_{i})):m_{A}(x_{i})\in [0,1]; \forall x_{i}\in X\right\}.
\end{eqnarray*}
FS $A$ on $X$ can be expressed by set of ordered pair as follows:
\begin{eqnarray*}
A=\left\{(x,m_{A}(x)):x\in X\right\}.
\end{eqnarray*}

SST developed by Molodtsov \cite{Molod} is a suitable tool for solving
uncertainties in non parametric situations and is a natural generalization of FS theory.
Since SST is a natural generalization of FS theory,
it has been applied in a wide range of fields ranging up to from mathematics to engineering from economics to optimization.
To deal with a collection of approximate description of objects, a generalized parametric gizmo is used known as SS.\\

In approximate description, there are two value sets which are called predicate and approximate.
Initially, the object description has an approximate by nature and so there is no require to present the concept of exact solution.
The SS theory is very convenient and simply effective in performance due to the nonentity of any limitations on the approximate descriptions.
With the aid of words and sentences, real sentences, real number, function, mapping and so on; any parameter can be operate that we desire.\\

\begin{Definition}\cite{Molod}
Consider $\mathcal{U}$, $\mathcal{P}$ as initial universe and parameters sets, respectively.
Take $\rho(U)$ as a power set of $\mathcal{U}$. Let $X\subset \mathcal{P}$.
Give the mapping $m:X \rightarrow \rho(\mathcal{U})$.
Therefore, $m_{X}$ is called a soft set(SS) on $\mathcal{U}$.
\end{Definition}

Choose set of $k$ objects and set of parameters as $\mathcal{U}=\{a_{1}, a_{2}, \ldots a_{k}\}$,  $\{A(1), A(2), \ldots, A(i)\}$, respectively.
Let $P \supseteq \{A(1)\cup A(2)\cup \cdots, A(i)\}$ and each parameter set $A(i)$ represent the
$i$th class of parameters and the elements of $A(i)$ represents a specific property set.
Assumed that the property sets can be shown as FSs.\\

\begin{Definition}\cite{Atanassov}
Let $\mathcal{U}$ be a universe. The set
\begin{eqnarray*}
\mathcal{A}=\{\langle x, m_{\mathcal{A}}(x), n_{\mathcal{A}}(x)\rangle: x\in \mathcal{U}\}
\end{eqnarray*}
is called an intuitionistic fuzzy set(IFS) $\mathcal{A}$ on $\mathcal{U}$,
where, $m_{\mathcal{A}}:\mathcal{U}\rightarrow [0,1]$ and $ n_{\mathcal{A}}:\mathcal{U}\rightarrow [0,1]$ such that $0\leq m_{\mathcal{A}}(x)+n_{\mathcal{A}}(x)\leq 1$
for any $x\in \mathcal{U}$.
\end{Definition}

The degree of indeterminacy $p_{\mathcal{A}}=1-m_{\mathcal{A}}(x)-n_{\mathcal{A}}(x)$.\\

\begin{Definition}\label{defFPSS}\cite{CagEn}
Let $\mathcal{U}$ be an initial universe, $\mathcal{P}$ be a set of all parameters and $X$ be a fuzzy set over $\mathcal{P}$.
The power set of $\mathcal{U}$ is denoted by $\rho(\mathcal{U})$.
If $ m_{X}:\mathcal{P}\rightarrow [0,1]$ and $ f_{X}:\mathcal{P}\rightarrow \rho(\mathcal{U})$ such that
$f_{X}=\emptyset$ if $m_{X}(x)=0$,
then, the set
\begin{eqnarray*}
\mathcal{U}_{X}=\{(m_{X}(x)/x, f_{X}(x)): x\in P\}
\end{eqnarray*}
is called fuzzy parameterized soft set(FPSS) on $\mathcal{U}$,
\end{Definition}

In Definition \ref{defFPSS}, $f_{X}$, $m_{X}$ called approximate function and
membership function of FPSS, respectively.

\begin{Definition}\label{defIFPSS}\cite{DeCa}
$X$ be an IFS over $\mathcal{P}$. An intuitionistic fuzzy parameterized sets(IFPS) $\mathcal{U}_{X}$ over $\mathcal{U}$
is defined as follows:
\begin{eqnarray*}
\mathcal{U}_{X}=\{(\langle x, \alpha_{X}(x), \beta_{X}(x)\rangle, f_{X}(x)): x\in \mathcal{P}\}
\end{eqnarray*}
where, $ \alpha_{X}:\mathcal{P}\rightarrow [0,1]$, $\beta_{X}:\mathcal{P}\rightarrow [0,1]$ and $f_{X}:\mathcal{P} \rightarrow \rho(\mathcal{U})$
with the property $f_{X}(x)=\emptyset$ if $\alpha_{X}(x)=0$ and $\beta_{X}(x)=1$.
\end{Definition}

In Definition \ref{defIFPSS}, the function $\alpha_{X}$ and $\beta_{X}$ called membership function and non-membership function
of IFPSS, respectively. The value $\alpha_{X}(x)$ and
$\beta_{X}(x)$ is degree of importance and unimportant of the parameter $x$.\\

Ordinary FPSS can be written as
\begin{eqnarray*}
\mathcal{U}_{X}=\{(\langle x, \alpha_{X}(x), 1-\alpha_{X}(x)\rangle, f_{X}(x)): x\in \mathcal{P}\}.
\end{eqnarray*}

\begin{Definition}\cite{Yager0, Yager1, YagerAbb}
An Pythagorean fuzzy set(PFS) $\varphi$ in $\mathcal{U}$ is given by
\begin{eqnarray*}
\varphi=\{\langle x, m_{\varphi}(x), n_{\varphi}(x)\rangle: x\in \mathcal{U}\},
\end{eqnarray*}
where $m_{\varphi}:\mathcal{U}\rightarrow [0,1]$ denotes the degree of membership and $n_{\varphi}:\mathcal{U}\rightarrow [0,1]$ denotes the degree of non-membership of the element $x\in \mathcal{U}$
to the set $\varphi$, respectively, with the condition that $0\leq (m_{\varphi}(x))^{2}+(n_{\varphi}(x))^{2}\leq 1$.
\end{Definition}

The degree of indeterminacy $\mathcal{I}_{\varphi}=\sqrt{1-(m_{\varphi}(x))^{2}-(n_{\varphi}(x))^{2}}$.\\

\begin{Definition}\cite{MajiBis1}
$X\subseteq \mathcal{P}$.
Then, $\varphi(X)$ is called Pythagorean Fuzzy Soft Set(PFSS)  on $\mathcal{U}$, if $\varphi(X)$ is mapping given by $\varphi(X):X\rightarrow \rho(\mathcal{U})$.
\end{Definition}

\begin{remark}
It is easy to check that PFSSs generalize both IFSs and SSs.
That is, all intuitionistic fuzzy degrees are part of the Pythagorean fuzzy degrees.
In actual decision-making problems, the PFSS characterizes a larger membership space than
the IFSS. Namely, the PFSS a higher capability
than the IFSS to model vagueness in real decision-making problems.\\
\end{remark}

Let $(\mathcal{L}, \leq_{\mathcal{L}})$ be a complete lattice, where $\mathcal{L}=\{(u,v): u,v \in [0,1], u^{2}+v^{2}<1\}$ and the corresponding partial order $\leq_{\mathcal{L}}$ is defined by
$(u,v)\leq_{\mathcal{L}} (i,j) \quad \Leftrightarrow \quad u\leq i \quad and \quad v \geq j$, for all $(u,v), (i,j)\in \mathcal{L}$. Any ordered pair $(u,v)\in \mathcal{L}$ is called Pythagorean fuzzy value(PFV) or
Pythagorean fuzzy number(PFN) \cite{garg}.\\

According to this new situation, we can represent PFS as follows:\\

For the $\mathcal{L}$-fuzzy set $\varphi:\mathcal{U} \rightarrow \mathcal{L}$, the PFS $\varphi(x)=\{(x, m_{\varphi}(x), n_{\varphi}(x)): x\in U\}$ can be identified as $\varphi(x)=(m_{\varphi}(x), n_{\varphi}(x))$  for all $x\in \mathcal{U}$.\\

\begin{figure} [htb]
\centering
  \includegraphics[width=100mm]{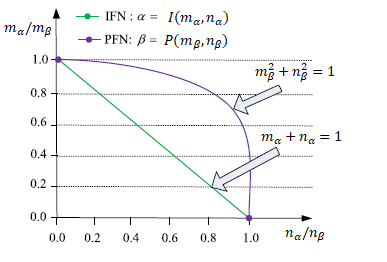}
  \caption{The PFNs and the IFNs }\label{Fig. 1}
\end{figure}

Let Pythagorean fuzzy numbers (PFNs) are denoted by $N=(m_{X}, n_{X})$ \cite{ZhangXu}.
Choose three PFNs $\theta = N(m, n), \theta_{1}=\langle m_{1}, n_{1}\rangle, \theta_{2}=\langle m_{2}, n_{2}\rangle$. We can give some basic operations as follows \cite{Yager0}, \cite{YagerAbb}:

\begin{itemize}
  \item $\bar{\theta}=\langle n, m \rangle$;
  \item $\theta_{1} \vee \theta_{2} = \langle \max\{m_{1}, m_{2}\} , \min\{n_{1}, n_{2}\}\rangle$;
  \item $\theta_{1} \wedge \theta_{2} = \langle \min\{m_{1}, m_{2}\} , \max\{n_{1}, n_{2}\}\rangle$;
  \item $\theta_{1} \oplus \theta_{2} = \langle \sqrt{m_{1}^{2}+m_{2}^{}2-m_{1}^{2}m_{2}^{2}} , n_{1}n_{2}\rangle$;
  \item $\theta_{1} \otimes \theta_{2} = \langle m_{1}m_{2}, \sqrt{n_{1}^{2}+n_{2}^{}2-n_{1}^{2}n_{2}^{2}} \rangle$;
  \item $\alpha.\theta = \langle \sqrt{1-(1-m^{2})^{\alpha}},n^{\alpha} \rangle$;
  \item $\theta ^{\alpha} = \langle m^{\alpha},  \sqrt{1-(1-n^{2})^{\alpha}}, \rangle$;
\end{itemize}

for $\alpha >0$.

\section{Pythagorean fuzzy parameterized soft set}

\begin{Definition}\label{def:0}
Let $\mathcal{U}$ be an initial universe, $\mathcal{P}$ be a set of all parameters and $X$ be Pythagorean fuzzy set on $\mathcal{P}$.
The set
\begin{eqnarray*}
\varphi_{X}=\bigg\{(\langle x,m_{X}(x), n_{X}(x)\rangle, f_{X}(x)): x\in \mathcal{P}\bigg\}
\end{eqnarray*}
is called an Pythagorean fuzzy parameterized soft set($\Phi-$soft set),
where $m_{X}:\mathcal{P}\rightarrow [0,1]$, $n_{X}:\mathcal{P}\rightarrow [0,1]$  and
$f_{X}:X\rightarrow \mathcal{L}$ is a Pythagorean fuzzy set such that
$m_{X}(x)=0$ and $n_{X}(x)=1$.
\end{Definition}

Here, the function $m_{X}$ and $n_{X}$ called membership function and
non-membership function of $\Phi-$soft set, respectively. The value
$m_{X}(x)$ and $n_{X}(x)$ is the degree of importance and unimportant of the parameter $x$. The elements of parameter $f_{X}$ are denoted by $(m_{f}, n_{f})$.\\

The set of all $\Phi-$soft set on $\mathcal{U}$ will be denoted by $\Phi(\mathcal{U})$.\\

The $\Phi-$soft set $\varphi_{X}$ on $\mathcal{U}$ can be represented by the set of ordered pairs,

\begin{eqnarray*}
\varphi_{X}=\left\{\left(\frac{x}{(m_{X}(x),n_{X}(x))}, f_{X}(x)\right): x\in P, f_{X}(x)\in \mathcal{L}, m_{X}(x),n_{X}(x)\in[0,1]\right\}.\\
\end{eqnarray*}

\begin{Example}\label{exp:1}
Let there be four patients in a clinic. Their symptoms are temperature $(s_{1})$, cough $(s_{2})$,
headache $(s_{3})$, chest problem $(s_{4})$, stomach problem $(s_{5})$ and myalgia $(s_{6})$.
The sets $\mathcal{P}=\{p_{1}, p_{2}, p_{3}, p_{4}\}$ and
$\mathcal{S}=\{s_{1}, s_{2}, s_{3}, s_{4}, s_{5}, s_{6}\}$ are patients and symptoms, respectively.
Let $X=\{s_{1}, s_{3}, s_{5}, s_{6}\}\subset \mathcal{S}$.
The values
\begin{eqnarray*}
f_{X}(s_{1})=(0.5, 0.4), \quad f_{X}(s_{3})=(0.7, 0.2), \quad f_{X}(s_{5})=(0.3, 0.6), \quad f_{X}(s_{6})=(0.6, 0.3)
\end{eqnarray*}
are state the diagnosis opinion of a physician.\\

All this information about patients can be represented in terms of the $\varphi_{X}$
as table in Table \ref{table:1}.

\begin{table}
 \caption{$\varphi_{X}$}\label{table:1}
  \vglue2mm
\centering
 {
\begin{tabular}{c  c  c  c  c }
\hline
P / X  & $s_{1}$ & $s_{3}$ & $s_{5}$ & $s_{6}$ \\
\hline
$p_{1}$  & (0.7, 0.7) & (0.6, 0.6) & (0.8, 0.6) & (0.4, 0.7)\\
$p_{2}$  & (0.5, 0.6) & (0.4, 0.5) & (0.8, 0.3) & (0.5, 0.6) \\
$p_{3}$ & (0.5, 0.4) & (0.9, 0.2) & (0.6, 0.4) & (0.6, 0.5) \\
$p_{4}$  & (0.7, 0.5) & (0.6, 0.2) & (0.5, 0.4) & (0.8, 0.4) \\
\hline
$f_{X}(s_{i})$ & (0.5, 0.4) &  (0.7, 0.2) & (0.3, 0.6) & (0.6, 0.3) \\
\hline
\end{tabular}}
\end{table}
\end{Example}

\begin{Definition}\label{def:1}
Let $\varphi_{X}, \varphi_{Y} \in \Phi(U)$. If the following conditions are hold, then
 $\varphi_{X}$ is $\Phi-$soft subset of $\varphi_{Y}$:
 \begin{itemize}
   \item [i.] $\varphi(X)(x)\subseteq_{\varphi} \varphi(Y)(x)$,
   \item [ii.] $m_{X}(x)\leq m_{Y}(x)$, $n_{X}(x)\geq n_{Y}(x)$
 \end{itemize}
for all $x\in P$.  $\Phi-$soft subset is denoted by $\varphi_{X} \hat{\subseteq} \varphi_{Y}$.
\end{Definition}

\begin{remark}
In this definition, it is not enough compare only $\varphi_{X}, \varphi_{Y}$ for $\varphi_{X} \hat{\subseteq} \varphi_{Y}$.
For $X,Y\subseteq P$ and $X \cap Y$, $f_{X}(x) \leq f_{Y}(x)$  may not be correct in all cases. Because the sets $X,Y$ are Pythagorean
fuzzy sets and their universes are different. Further, since $X \cap Y$, the sets $X,Y$ have not common parameter $x$.
So it has to be $\varphi(X)\subseteq_{\varphi} \varphi(X)$. That is,
if the inclusion relation $f_{X}(x) \leq f_{Y}(x)$ is holds on Definition \ref{def:1}, then
$\varphi_{X} \hat{\subseteq} \varphi_{Y}$ may not be true for every element, in contrast to the definition of classical subset.
Therefore, conditions (i) and (ii) are mandatory for definition of $\Phi-$soft subset.
\end{remark}

\begin{Definition}\label{def:2}
Let $\varphi_{X}, \varphi_{Y} \in \Phi(U)$. Then, $\varphi_{X}=\varphi_{Y}$ if
$X=Y$, $\varphi(X)=\varphi(Y)$ and $f_{X}=f_{Y}$.
\end{Definition}

From Definition \ref{def:1} and Definition \ref{def:2}, the following Proposition can easily be proved:

\begin{Proposition}
Let $\varphi_{X}, \varphi_{Y}, \varphi_{Z} \in \Phi(U)$. Then, the following conditions are hold:
\begin{itemize}
\item[i.] $\varphi_{X} \hat{\subseteq} \varphi_{Y}$ and $\varphi_{Y} \hat{\subseteq} \varphi_{Z} \Leftrightarrow \varphi_{X} \hat{\subseteq} \varphi_{Z}$
\item[ii.] $\varphi_{X} \hat{\subseteq} \varphi_{Y}$ and $\varphi_{Y} \hat{\subseteq} \varphi_{X} \Leftrightarrow \varphi_{X}=\varphi_{Y}$
\item[iii.] $\varphi_{X}=\varphi_{Y}$ and $\varphi_{Z}=\varphi_{X} \Leftrightarrow \varphi_{X}=\varphi_{Z}$
\end{itemize}
\end{Proposition}

\section{Operations of $\Phi(U)$}

\begin{Definition}\label{def:4}
Suppose that $X,Y \subseteq P$ and $Z=X \cup Y$.
Let $\varphi_{X}, \varphi_{Y} \in \Phi(U)$. The extended union of $\varphi_{X}, \varphi_{Y}$ is defined as
$\varphi_{Z}=\varphi_{X} \hat{\cup}_{E} ~\varphi_{Y}$ such that for all $z\in Z$,

\begin{eqnarray*}
\varphi(X)\cup_{E} ~\varphi(Y)=\varphi(Z);
\end{eqnarray*}
\begin{eqnarray*}
m_{f_{Z}}(z)= \left\{ \begin{array}{ccl}
m_{f_{X}}(z)&, & \quad z\in X - Y, \\
m_{f_{Y}}(z)&, & \quad z\in Y - X, \\
\max\{m_{f_{X}}(z), m_{f_{Y}}(z)\}&, & z\in X\cap Y;
\end{array} \right.
\end{eqnarray*}

\begin{eqnarray*}
n_{f_{Z}}(z)= \left\{ \begin{array}{ccl}
n_{f_{X}}(z)&, & \quad z\in X - Y, \\
n_{f_{Y}}(z)&, & \quad z\in Y - X, \\
\min\{n_{f_{X}}(z), n_{f_{Y}}(z)\}&, & z\in X\cap Y;
\end{array} \right.
\end{eqnarray*}
\end{Definition}

\begin{Definition}\label{def:5}
Suppose that $X,Y \subseteq P$ and $Z=X \cup Y$.
Let $\varphi_{X}, \varphi_{Y} \in \Phi(U)$. The extended intersection of $\varphi_{X}, \varphi_{Y}$ is defined as
$\varphi_{Z}=\varphi_{X} \hat{\cap}_{E} ~\varphi_{Y}$ such that for all $z\in Z$,

\begin{eqnarray*}
\varphi(X)\cap_{E} ~\varphi(Y)=\varphi(Z);
\end{eqnarray*}
\begin{eqnarray*}
m_{f_{Z}}(z)= \left\{ \begin{array}{ccl}
m_{f_{X}}(z)&, & \quad z\in X - Y, \\
m_{f_{Y}}(z)&, & \quad z\in Y - X, \\
\min\{m_{f_{X}}(z), m_{f_{Y}}(z)\}&, & z\in X\cap Y;
\end{array} \right.
\end{eqnarray*}

\begin{eqnarray*}
n_{f_{Z}}(z)= \left\{ \begin{array}{ccl}
n_{f_{X}}(z)&, & \quad z\in X - Y, \\
n_{f_{Y}}(z)&, & \quad z\in Y - X, \\
\max\{n_{f_{X}}(z), n_{f_{Y}}(z)\}&, & z\in X\cap Y;
\end{array} \right.
\end{eqnarray*}
\end{Definition}

\begin{remark}
In Definitions \ref{def:4} and \ref{def:5}, since $\varphi(X), \varphi(Y) $ are Pythagorean fuzzy sets, the operations t-norm and t-conorm are not used.
That is, t-norm and t-conorm are binary functions on the $[0,1]$ and so the notations $\varphi(X) \circ \varphi(Y) $, $\varphi(X) \ast \varphi(Y) $
are not used, where the operations $\circ$, $\ast$ are represent t-norm and t-conorm, respectively. Similarly, if $S$ is a Pyhtagorean fuzzy set and
 $f(X), f(Y)\in S$,
then the notations $f(X) \circ f(Y)$, $f(X)\ast f(Y)$ are also  improperly used.
\end{remark}

Now, we give restricted union and restricted intersection of $\varphi_{X}, \varphi_{Y} $.
These operations are denoted by $\hat{\cup}_{R}$ and  $\hat{\cap}_{R}$, respectively.

\begin{Definition}\label{def:40}
Choose $X,Y \subseteq P$ and $Z(R)=X \cap Y\neq \emptyset$. Let $\varphi_{X}, \varphi_{Y} \in \Phi(U)$.
The restricted union of $\varphi_{X}, \varphi_{Y}$ is defined as
$\varphi_{Z(R)}=\varphi_{X} \hat{\cup}_{R} ~\varphi_{Y}$ such that
\begin{eqnarray*}
\varphi(Z(R))=\varphi(X)\cup_{R} ~\varphi(Y);
\end{eqnarray*}
and for all $z\in Z$,
\begin{eqnarray*}
m_{f_{Z(R)}}(z)=\max\{m_{f_{X}}(z), m_{f_{Y}}(z)\}; \quad n_{f_{Z(R)}}(z)=\min\{n_{f_{X}}(z), n_{f_{Y}}(z)\}.
\end{eqnarray*}
\end{Definition}

\begin{Definition}\label{def:50}
Choose $X,Y \subseteq P$ and $T(R)=X \cap Y\neq \emptyset$. Let $ \varphi_{X}, \varphi_{Y}\in \Phi(U)$.
The restricted intersection of $\varphi_{X}, \varphi_{Y}$ is defined as
$\varphi_{T(R)}=\varphi_{X} \hat{\cap}_{R} ~\varphi_{Y}$ such that
\begin{eqnarray*}
\varphi(T(R))=\varphi(X)\cup_{R} ~\varphi(Y);
\end{eqnarray*}
and for all $t\in T$,
\begin{eqnarray*}
m_{f_{T(R)}}(t)=\min\{m_{f_{X}}(t), m_{f_{Y}}(t)\}; \quad n_{f_{T(R)}}(t)=\max\{n_{f_{X}}(t), n_{f_{Y}}(t)\}.
\end{eqnarray*}
\end{Definition}

\begin{Example}\label{exp:2}
Take the sets $\mathcal{P}, \mathcal{S}, X$ as in Example \ref{exp:1}. Choose the set $Y=\{s_{2}, s_{3}, s_{5}, s_{6}\}\subset \mathcal{S}$.
For the $\varphi_{X}, \varphi_{Y}$, consider the extended union and the extended intersection
\begin{eqnarray*}
\varphi_{Z}^{\cup}=\varphi_{X} \hat{\cup}_{E} ~\varphi_{Y} \quad  and \quad \varphi_{T}^{\cap}=\varphi_{X} \hat{\cap}_{E} ~\varphi_{Y},
\end{eqnarray*}
where $\varphi_{Z}^{\cup}=(Z, X\cup Y, f_{Z})$, ~$\varphi_{T}^{\cap}=(T, X\cup Y, f_{T})$.\\

We can calculate the extended union as follows (Table \ref{table:3}):
\begin{eqnarray*}
&& Z(s_{1})=\left\{\frac{p_{1}}{(0.7,0.7)}, \frac{p_{2}}{(0.5,0.6)}, \frac{p_{3}}{(0.5,0.4)}, \frac{p_{4}}{(0.7,0.5)}\right\},\\
&& Z(s_{2})=\left\{\frac{p_{1}}{(0.6,0.6)}, \frac{p_{2}}{(0.1,0.7)}, \frac{p_{3}}{(0.3,0.4)}, \frac{p_{4}}{(0.5,0.4)}\right\},\\
&& Z(s_{3})=\left\{\frac{p_{1}}{(0.6,0.2)}, \frac{p_{2}}{(0.4,0.5)}, \frac{p_{3}}{(0.9,0.2)}, \frac{p_{4}}{(0.6,0.2)}\right\},\\
&& Z(s_{5})=\left\{\frac{p_{1}}{(0.8,0.4)}, \frac{p_{2}}{(0.8,0.1)}, \frac{p_{3}}{(0.6,0.4)}, \frac{p_{4}}{(0.6,0.4)}\right\},\\
&& Z(s_{6})=\left\{\frac{p_{1}}{(0.4,0.5)}, \frac{p_{2}}{(0.5,0.5)}, \frac{p_{3}}{(0.6,0.2)}, \frac{p_{4}}{(0.8,0.5)}\right\}
\end{eqnarray*}
and $f_{Z}(s_{1})=(0.5, 0.4)$, $f_{Z}(s_{2})=(0.1, 0.6)$, $f_{Z}(s_{3})=(0.7, 0.2)$, $f_{Z}(s_{5})=(0.4, 0.5)$, $f_{Z}(s_{6})=(0.6, 0.3)$.\\

Now, we calculate the extended intersection (Table \ref{table:4}):
\begin{eqnarray*}
&& T(s_{1})=\left\{\frac{p_{1}}{(0.7,0.7)}, \frac{p_{2}}{(0.5,0.6)}, \frac{p_{3}}{(0.5,0.4)}, \frac{p_{4}}{(0.7,0.5)}\right\},\\
&& T(s_{2})=\left\{\frac{p_{1}}{(0.6,0.6)}, \frac{p_{2}}{(0.1,0.7)}, \frac{p_{3}}{(0.3,0.4)}, \frac{p_{4}}{(0.5,0.4)}\right\},\\
&& T(s_{3})=\left\{\frac{p_{1}}{(0.4,0.6)}, \frac{p_{2}}{(0.3,0.5)}, \frac{p_{3}}{(0.7,0.4)}, \frac{p_{4}}{(0.5,0.2)}\right\},\\
&& T(s_{5})=\left\{\frac{p_{1}}{(0.6,0.6)}, \frac{p_{2}}{(0.5,0.3)}, \frac{p_{3}}{(0.2,0.5)}, \frac{p_{4}}{(0.5,0.4)}\right\},\\
&& T(s_{6})=\left\{\frac{p_{1}}{(0.1,0.7)}, \frac{p_{2}}{(0.2,0.6)}, \frac{p_{3}}{(0.4,0.5)}, \frac{p_{4}}{(0.5,0.5)}\right\}
\end{eqnarray*}
and $f_{T}(s_{1})=(0.5, 0.4)$, $f_{T}(s_{2})=(0.1, 0.6)$, $f_{T}(s_{3})=(0.7, 0.2)$, $f_{T}(s_{5})=(0.3, 0.6)$, $f_{T}(s_{6})=(0.6, 0.3)$.\\
\end{Example}

\begin{table}
 \caption{$\varphi_{Y}$}\label{table:2}
  \vglue2mm
\centering
 {
\begin{tabular}{c  c  c  c  c }
\hline
P / Y  & $s_{2}$ & $s_{3}$ & $s_{5}$ & $s_{6}$ \\
\hline
$p_{1}$  & (0.6, 0.6) & (0.4, 0.2) & (0.6, 0.4) & (0.1, 0.5)\\
$p_{2}$  & (0.1, 0.7) & (0.3, 0.5) & (0.5, 0.1) & (0.2, 0.5) \\
$p_{3}$ & (0.3, 0.4) & (0.7, 0.4) & (0.2, 0.5) & (0.4, 0.2) \\
$p_{4}$  & (0.5, 0.4) & (0.5, 0.2) & (0.6, 0.4) & (0.5, 0.5) \\
\hline
$f_{Y}(s_{i})$ & (0.1, 0.6) &  (0.7, 0.2) & (0.4, 0.5) & (0.6, 0.3) \\
\hline
\end{tabular}}
\end{table}

Using the information in the Example \ref{exp:1} and Example \ref{exp:2} examples,
we can obtain the restricted union and restricted intersection, as follows (Table \ref{table:5}, Table \ref{table:6}):

\begin{eqnarray*}
&& Z(R)(s_{3})=\left\{\frac{p_{1}}{(0.6,0.2)}, \frac{p_{2}}{(0.4,0.5)}, \frac{p_{3}}{(0.9,0.2)}, \frac{p_{4}}{(0.6,0.2)}\right\},\\
&& Z(R)(s_{5})=\left\{\frac{p_{1}}{(0.8,0.4)}, \frac{p_{2}}{(0.8,0.1)}, \frac{p_{3}}{(0.6,0.4)}, \frac{p_{4}}{(0.6,0.4)}\right\},\\
&& Z(R)(s_{6})=\left\{\frac{p_{1}}{(0.4,0.5)}, \frac{p_{2}}{(0.5,0.5)}, \frac{p_{3}}{(0.6,0.2)}, \frac{p_{4}}{(0.8,0.5)}\right\}
\end{eqnarray*}
and $f_{Z}(s_{3})=(0.7, 0.2)$, $f_{Z}(s_{5})=(0.4, 0.5)$, $f_{Z}(s_{6})=(0.6, 0.3)$.

\begin{eqnarray*}
&& T(R)(s_{3})=\left\{\frac{p_{1}}{(0.4,0.6)}, \frac{p_{2}}{(0.3,0.5)}, \frac{p_{3}}{(0.7,0.4)}, \frac{p_{4}}{(0.5,0.2)}\right\},\\
&& T(R)(s_{5})=\left\{\frac{p_{1}}{(0.6,0.6)}, \frac{p_{2}}{(0.5,0.3)}, \frac{p_{3}}{(0.2,0.5)}, \frac{p_{4}}{(0.5,0.4)}\right\},\\
&& T(R)(s_{6})=\left\{\frac{p_{1}}{(0.1,0.7)}, \frac{p_{2}}{(0.2,0.6)}, \frac{p_{3}}{(0.4,0.5)}, \frac{p_{4}}{(0.5,0.5)}\right\}
\end{eqnarray*}
and $f_{T}(s_{3})=(0.7, 0.2)$, $f_{T}(s_{5})=(0.3, 0.6)$, $f_{T}(s_{6})=(0.6, 0.3)$.\\

\begin{table}
 \caption{Extended Union}\label{table:3}
  \vglue2mm
\centering
 {
\begin{tabular}{c  c  c  c  c c}
\hline
P / Z  & $s_{1}$ & $s_{2}$ & $s_{3}$ & $s_{5}$ & $s_{6}$ \\
\hline
$p_{1}$  & (0.7, 0.7) & (0.6, 0.6) & (0.6, 0.2) & (0.8, 0.4) & (0.4, 0.5)\\
$p_{2}$  & (0.5, 0.6) & (0.1, 0.7) & (0.4, 0.5) & (0.8, 0.1) & (0.5, 0.5) \\
$p_{3}$ & (0.5, 0.4) & (0.3, 0.4) & (0.9, 0.4) & (0.6, 0.4) & (0.6, 0.2) \\
$p_{4}$  & (0.7, 0.5) & (0.5, 0.4) & (0.6, 0.2) & (0.6, 0.4) & (0.8, 0.5) \\
\hline
$f_{Z}(s_{i})$ & (0.5, 0.4) & (0.1, 0.6) &  (0.7, 0.2) & (0.4, 0.5) & (0.6, 0.3) \\
\hline
\end{tabular}}
\end{table}

\begin{table}
 \caption{Extended Intersection}\label{table:4}
  \vglue2mm
\centering
 {
\begin{tabular}{c  c  c  c  c c}
\hline
P / T  & $s_{1}$ & $s_{2}$ & $s_{3}$ & $s_{5}$ & $s_{6}$ \\
\hline
$p_{1}$  & (0.7, 0.7) & (0.6, 0.6) & (0.4, 0.6) & (0.6, 0.6) & (0.1, 0.7)\\
$p_{2}$  & (0.5, 0.6) & (0.1, 0.7) & (0.3, 0.5) & (0.5, 0.3) & (0.2, 0.6) \\
$p_{3}$ & (0.5, 0.4) & (0.3, 0.4) & (0.7, 0.4) & (0.2, 0.5) & (0.4, 0.5) \\
$p_{4}$  & (0.7, 0.5) & (0.5, 0.4) & (0.5, 0.2) & (0.5, 0.4) & (0.5, 0.5) \\
\hline
$f_{T}(s_{i})$ & (0.5, 0.4) & (0.1, 0.6) &  (0.7, 0.2) & (0.3, 0.6) & (0.6, 0.3) \\
\hline
\end{tabular}}
\end{table}

\begin{table}
 \caption{Restricted Union}\label{table:5}
  \vglue2mm
\centering
 {
\begin{tabular}{c  c  c  c  }
\hline
P / Z(R)  & $s_{3}$ & $s_{5}$ & $s_{6}$ \\
\hline
$p_{1}$  & (0.6, 0.2) & (0.8, 0.4) & (0.4, 0.5) \\
$p_{2}$  & (0.4, 0.5) & (0.8, 0.1) & (0.5, 0.5) \\
$p_{3}$ & (0.9, 0.2) & (0.6, 0.4) & (0.6, 0.2)  \\
$p_{4}$  & (0.6, 0.2) & (0.4, 0.5) & (0.6, 0.3)  \\
\hline
$f_{Z(R)}(s_{i})$ & (0.7, 0.2) & (0.4, 0.5) &  (0.6, 0.3)  \\
\hline
\end{tabular}}
\end{table}

\begin{table}
 \caption{Restricted Intersection}\label{table:6}
  \vglue2mm
\centering
 {
\begin{tabular}{c  c  c  c  }
\hline
P / T(R)  & $s_{3}$ & $s_{5}$ & $s_{6}$ \\
\hline
$p_{1}$  & (0.4, 0.6) & (0.6, 0.6) & (0.1, 0.7) \\
$p_{2}$  & (0.3, 0.5) & (0.5, 0.3) & (0.2, 0.6) \\
$p_{3}$ & (0.7, 0.4) & (0.2, 0.5) & (0.4, 0.5)  \\
$p_{4}$  & (0.5, 0.2) & (0.5, 0.4) & (0.5, 0.5)  \\
\hline
$f_{T(R)}(s_{i})$ & (0.7, 0.2) & (0.3, 0.6) &  (0.6, 0.3)  \\
\hline
\end{tabular}}
\end{table}

Let $\varphi(X) \in \Phi(U)$ and $a,b \in [0,1]$, such that $a^{2}+b^{2}\leq 1$. Then, $\varphi_{X}$ is called
a $(a,b)-$constant PFSS, denoted by $X^{(a,b)}$, if $m_{\varphi_{X}(x)}=\tilde{a}$ and
$n_{\varphi_{X}(x)}=\tilde{b}$ for all $x\in X$.\\

In this text, $X^{(0,1)}$ and $X^{(1,0)}$ are represented the relative null PFSS and relative whole PFSS
with respect to the parameter set $P$, respectively.\\

\begin{Definition}\label{def:60}
Let $\varphi_{X} \in \Phi(U)$.
\begin{itemize}
\item[i.] $\varphi_{X}$ is called a relative null $\Phi-$soft set with respect to $X$,  if $\varphi(X)=X^{(0,1)}$ and $m_{f_{X}}(x)=0, n_{f_{X}}(x)=1$  for all $x\in X$.
\item[ii.] $\varphi_{X}$ is called a relative whole $\Phi-$soft set, with respect to $X$, if $\varphi(X)=X^{(1,0)}$ and $m_{f_{X}}(x)=1, n_{f_{X}}(x)=0$  for all $x\in X$.
\end{itemize}
\end{Definition}

The relative null $\Phi-$soft set and relative whole $\Phi-$soft set are denoted by $\varphi_{\emptyset}$ and $\varphi_{\hat{X}}$, respectively. Furthermore, if we choose $X=P$, then
$\varphi_{\emptyset}$ and $\varphi_{\hat{X}}$ are called the null $\Phi-$soft set and whole $\Phi-$soft set, respectively.\\

From the definitions \ref{def:0}, \ref{def:4}, \ref{def:5}, \ref{def:40}, \ref{def:50} and \ref{def:60}, we may give the Proposition \ref{prop:1}:

\begin{Proposition}\label{prop:1}
Let $\varphi_{X}\in \Phi(U)$. Then, the following conditions are hold:
\begin{itemize}
\item[i.] $\varphi_{X} \hat{\cup}_{E} ~\varphi_{X}=\varphi_{X} \hat{\cup}_{R} ~\varphi_{X}=\varphi_{X}$,
\item[ii.] $\varphi_{X} \hat{\cap}_{E} ~\varphi_{X}=\varphi_{X} \hat{\cap}_{R} ~\varphi_{X}=\varphi_{X}$,
\item[iii.] $\varphi_{X} \hat{\cup}_{E} ~\varphi_{\emptyset}=\varphi_{X} \hat{\cup}_{R} ~\varphi_{\emptyset}=\varphi_{X}$,
\item[iv.] $\varphi_{X} \hat{\cap}_{E} ~\varphi_{\emptyset}=\varphi_{X} \hat{\cap}_{R} ~\varphi_{\emptyset}=\varphi_{\emptyset}$
\item[v.] $\varphi_{X} \hat{\cup}_{E} ~\varphi_{\hat{X}}=\varphi_{X} \hat{\cup}_{R} ~\varphi_{\hat{X}}=\varphi_{\hat{X}}$,
\item[vi.] $\varphi_{X} \hat{\cap}_{E} ~\varphi_{\hat{X}}=\varphi_{X} \hat{\cap}_{R} ~\varphi_{\hat{X}}=\varphi_{X}$.
\end{itemize}
\end{Proposition}


\section{$\Phi-$Soft Decision-Making Method}

\subsection{Score and Accuracy Functions}

For PFNs, the mapping $\mathcal{SF}:\mathcal{L}\rightarrow [-1,1]$ is called \emph{score function}, if
\begin{eqnarray*}
\mathcal{SF}_{N}=m_{N}^{2}-n_{N}^{2}
\end{eqnarray*}
for all $N=(m_{N}, n_{N})\in \mathcal{L}$ \cite{ZhangXu}, \cite{Agarwaletal}. \\

For any two PFNs $N,M$, $N \prec M$ if $\mathcal{SF}(N)< \mathcal{SF}(M)$; $N \succ M$ if $\mathcal{SF}(N)> \mathcal{SF}(M)$;
$N \sim M$ if $\mathcal{SF}(N) = \mathcal{SF}(M)$.\\

As can be seen from the definition of $\mathcal{SF}$,
the score the larger the score $s_{N}$, the greater the PFN $N$ \cite{garg}.
It should be noted that $\mathcal{SF}$ cannot differentiate some evidently distinct
PFNs which have the same score.
We can give examples to explain this situation:
Take two PFNs $N,M$ as $N=(0.481, 0.402)$ and $M=(0.527, 0.456)$.
Then, $SF_{N}=0.0697$ and $SF_{M}=0.0697$. Again, for $N=(0.123, 0.123)$ and $M=(0.456, 0.456)$
we can write $N \sim M$.
Therefore, if only the scoring function is used for comparison, it is not possible to make a comparison between these numbers. \\

To overcome this problem, we can define a new function \cite{Agarwaletal}, as follows:\\

The mapping $\mathcal{AF}:\mathcal{L}\rightarrow [0,1]$ is called \emph{accuracy function}, if
\begin{eqnarray*}
\mathcal{AF}_{N}=m_{N}^{2}+n_{N}^{2}
\end{eqnarray*}
for all $N=(m_{N}, n_{N})\in \mathcal{L}$ \cite{PengYang}. \\

Using the $\mathcal{SF}$ and $\mathcal{AF}$, for comparing PFVs, the following method is presented by Agarwal et al. \cite{Agarwaletal}.\\

For any two PFNs $N,M$, if $\mathcal{SF}(N) = \mathcal{SF}(M)$, then (i) $N>M$ if $\mathcal{AF}(N)> \mathcal{AF}(M)$ (ii) $N<M$ if $\mathcal{AF}(N)< \mathcal{AF}(M)$and
(iii) $N \sim M$ if $\mathcal{AF}(N)= \mathcal{AF}(M)$. For a binary relation $\leq_{(\mathcal{SF}, \mathcal{AH})}\in \mathcal{L}$ and $N,M \in \mathcal{L}$, it can be written as

\begin{eqnarray*}
N \leq_{(\mathcal{SF}, \mathcal{AH})} M \quad \Leftrightarrow \quad (\mathcal{SF}_{N}< \mathcal{SF}_{M})\vee (\mathcal{SF}_{N}= \mathcal{SF}_{M} \wedge \mathcal{AF}_{N} \leq \mathcal{AF}_{M}).
\end{eqnarray*}

From \cite{ChenTan} and \cite{Fengetal}, we can give the new definition:\\

\begin{Definition}
The mapping $\mathcal{ES}: \mathcal{L}\rightarrow [0,1]$ is called expectation score function such that for all $N=(m_{N}, n_{N})\in \mathcal{L}$
\begin{eqnarray*}
\mathcal{ES}_{N}=\frac{m_{N}^{2}-n_{N}^{2}+1}{2}.
\end{eqnarray*}
\end{Definition}

\begin{remark}
In this definition, if we take $(m_{N}^{2})^{*}=1-n_{N}^{2}$, then,
\begin{eqnarray*}
\mathcal{ES}_{N}=\frac{m_{N}^{2}+(1-n_{N}^{2})}{2}=\frac{m_{N}^{2}+(m_{N}^{2})^{*}}{2}.
\end{eqnarray*}
Therefore, the PFV $N=(m_{N}, n_{N})$ specifies an interval $[m_{N}, (m_{N})^{*}]$, in which
lies the accurate value of membership grade. Hence, the uncertain value of the membership grade can be considered
as a random variable $x_{N}$. Assume that $x_{N}$ be uniformly distributed on the $[m_{N}, (m_{N})^{*}]$. Then,
the value of expectation score function $ES$ is obtained as $ES_{N}=exc(x_{N})$. Further, the function $ES$ is bounded between 0 and 1.\\
\end{remark}

\begin{Proposition}
Let $\mathcal{ES}: \mathcal{L}\rightarrow [0,1]$ and $N=(m_{N}, n_{N}) \in \mathcal{L}$. Then, we have
\begin{itemize}
\item [i.] $\mathcal{ES}(0,1)=0$ and  $\mathcal{ES}(1,0)=1$,
\item [ii.] $\mathcal{ES}(m_{N}, n_{N})$ is increasing with respect to $m_{N}$,
\item [iii.] $\mathcal{ES}(m_{N}, n_{N})$ is decreasing with respect to $n_{N}$.
\end{itemize}
\end{Proposition}

From the Definition $\mathcal{ES}$ and the remark, this Proposition can be easily proved.\\

We can write the following definitions and theorem with the same idea in Feng et al. \cite{Fengetal}.\\

\begin{Definition}\label{def:80}
For two PFNs $N, M \in \mathcal{L}$ and the relation $\leq_{(m,\mathcal{ES})}$ on $\mathcal{L}$, we have
\begin{eqnarray*}
N \leq_{(m,\mathcal{ES})} M \quad \Leftrightarrow \quad (m_{N}< m_{M})\vee (m_{N}= m_{M} \wedge \mathcal{ES}_{N} \leq \mathcal{ES}_{M}).
\end{eqnarray*}
\end{Definition}

\begin{Theorem}
The relation $\leq_{(m,\mathcal{ES})}$ is a partial order on $\mathcal{L}$.
\end{Theorem}

By replacing the approval rates in Definition \ref{def:80}, the other relation $\leq_{(\mathcal{ES},m)}$ is written as follows:\\

\begin{Definition}
For two PFNs $N, M \in \mathcal{L}$ and the relation $\leq_{(m,\mathcal{ES})}$ on $\mathcal{L}$,
\begin{eqnarray*}
N \leq_{(m,\mathcal{ES})} M \quad \Leftrightarrow \quad (\mathcal{ES}_{N}< \mathcal{ES}_{M})\vee (\mathcal{ES}_{N}= \mathcal{ES}_{M} \wedge m_{N} \leq m_{M}).
\end{eqnarray*}
\end{Definition}

The relationship with each other of $\leq_{(\mathcal{SF},\mathcal{AF})}$ and $\leq_{(\mathcal{ES},m)}$ can be given as follows:

\begin{Proposition}
Let $N$ and $M$ be PFNs in $\mathcal{L}$. Then, $N \leq_{(\mathcal{SF},\mathcal{AF})} M$ iff $N \leq_{(\mathcal{ES},m)} M$.
\end{Proposition}

\begin{proof}
Assume that $N \leq_{(\mathcal{SF},\mathcal{AF})} M$. Then, two cases are valid:\\

i. For $\mathcal{SF}_{N} < \mathcal{SF}_{M}$, $\mathcal{ES}_{N} < \mathcal{ES}_{M}$ becomes. Therefore, $N \leq_{(\mathcal{ES},m)} M$.\\

ii.  $\mathcal{ES}_{N} = \mathcal{ES}_{M}$ becomes, for $\mathcal{SF}_{N} = \mathcal{SF}_{M}$ and $\mathcal{AF}_{N} \leq \mathcal{AF}_{M}$. Moreover,
$\mathcal{AF}_{M}-\mathcal{AF}_{N} \geq 0$. Further, $\mathcal{SF}_{N} = \mathcal{SF}_{M} \Rightarrow n_{M}^{2}-n_{N}^{2}=m_{M}^{2}-m_{N}^{2}$. It follows that $m_{M}^{2}-m_{N}^{2}=(\mathcal{AF}_{N} - \mathcal{AF}_{M})/2 \geq 0$. In that case, $N \leq_{(\mathcal{ES},m)} M$.\\

Now, we take $N \leq_{(\mathcal{ES},m)} M$. In this also situation two cases are hold:\\

i. For $\mathcal{ES}_{N} < \mathcal{ES}_{M}$, $\mathcal{SF}_{N} < \mathcal{SF}_{M}$ and therefore $N \leq_{(\mathcal{SF},\mathcal{AF})} M$.\\

ii. $\mathcal{SF}_{N} = \mathcal{SF}_{M}$ becomes, for $\mathcal{ES}_{N} = \mathcal{ES}_{M}$ and $m_{N} \leq m_{M}$. Then,
$n_{M}^{2}-n_{N}^{2}=m_{M}^{2}-m_{N}^{2}$. It follows that $\mathcal{AF}_{M}-\mathcal{AF}_{N}=(m_{N}^{2} - m_{M}^{2})+(n_{N}^{2} - n_{M}^{2})=2(m_{N}^{2} - m_{M}^{2})\geq 0$. In that case, $N \leq_{(\mathcal{SF},\mathcal{AF})} M$.
\end{proof}

If we consider the debates of \cite{Xu} and \cite{Fengetal}, then the following proposition is given:\\

\begin{Proposition}
Take two PFVs $N=(m_{N}, n_{N})$ and $M=(m_{M}, n_{M})$ in $\mathcal{L}$.
Then, the following conditions are equivalent:
\begin{itemize}
\item [i.] $\mathcal{SF}_{N}=\mathcal{SF}_{M} \wedge \mathcal{AF}_{N} \leq \mathcal{AF}_{M}$,
\item [ii.] $\mathcal{ES}_{N}=\mathcal{ES}_{M} \wedge m_{N} \leq m_{M}$,
\item [ii.] $\mathcal{ES}_{N}=\mathcal{ES}_{M} \wedge n_{N} \leq n_{M}$,
\item [iv.] $\mathcal{SF}_{N}=\mathcal{SF}_{M} \wedge m_{N} \leq m_{M}$,
\item [v.] $\mathcal{SF}_{N}=\mathcal{SF}_{M} \wedge n_{N} \leq n_{M}$,
\end{itemize}
\end{Proposition}

\begin{proof}
Consider the condition iv and choose two PFVs $N=(m_{N}, n_{N}), M=(m_{M}, n_{M})\in \mathcal{L}$.
Since $\mathcal{SF}_{N}=\mathcal{SF}_{M}$, we have $m_{N}^{2} - m_{M}^{2}= n_{N}^{2} - n_{M}^{2}$.
It follows that $m_{N} \leq m_{M}$ iff $n_{N} \leq n_{M}$. Therefore, the condition iv equivalent
to the condition v.\\

Similarly, other equivalents will be proved, so we omit it.
\end{proof}

\begin{Definition}\cite{ZhangXu}\label{def:70}
Let $N=(m_{N}, n_{N}), M=(m_{M}, n_{M})\in \mathcal{L}$ be two PFVs. Then, for $\alpha>0$, we have the following operations:\\
\begin{itemize}
\item [i.] $M +_{P} N = \left(\sqrt{m_{M}^{2}+m_{N}^{2}-m_{M}^{2}m_{N}^{2}}, n_{M}^{2}n_{N}^{2}\right)$,
\item [ii.] $M \times_{P} N= \left(m_{M}^{2}m_{N}^{2}, \sqrt{n_{M}^{2}+n_{N}^{2}-n_{M}^{2}n_{N}^{2}} \right)$,
\item [ii.] $\alpha N= \left(\sqrt{1-(1-m_{N}^{2})^{\alpha}}, (n_{N})^{\alpha}\right)$,
\item [iv.] $N^{\alpha} = \left((m_{N})^{\alpha}, \sqrt{1-(1-n_{N}^{2})^{\alpha}}, \right)$.
\end{itemize}
\end{Definition}

\begin{Theorem}\cite{ZhangXu}
For $N=(m_{N}, n_{N}), M=(m_{M}, n_{M})\in \mathcal{L}$ and $\alpha, \alpha_{1}, \alpha_{2}>0$,\\
\begin{itemize}
\item [i.] $M +_{P} N = N +_{P} M$,
\item [ii.] $M \times_{P} N= N \times_{P} M$
\item [iii.] $\alpha(M+N)=\alpha M +_{P} \alpha N$,
\item [iv.] $\alpha_{1}M +_{P} \alpha_{2}M=(\alpha_{1}+ \alpha_{2})M$,
\item [v.] $(M \times_{P} N)^{\alpha}=M^{\alpha}\times_{P} N^{\alpha}$,
\item [vi.] $M^{\alpha_{1}} \times_{P} M^{\alpha_{2}}=M^{(\alpha_{1}+\alpha_{2})}$.
\end{itemize}
\end{Theorem}

\begin{Proposition}
For $M=(m_{M}, n_{M}), N=(m_{N}, n_{N}), K=(m_{K}, n_{K})\in \mathcal{L}$, then we have,
\begin{eqnarray*}
N \leq_{(m,\mathcal{ES})} K \Rightarrow M +_{P} N \leq_{(m,\mathcal{ES})} M +_{P} K.
\end{eqnarray*}
\end{Proposition}

\begin{proof}
From Definition \ref{def:70}, we can write, $M +_{P} N = \left(\sqrt{m_{M}^{2}+m_{N}^{2}-m_{M}^{2}m_{N}^{2}}, n_{M}^{2}n_{N}^{2}\right)$
and $M +_{P} K = \left(\sqrt{m_{M}^{2}+m_{K}^{2}-m_{M}^{2}m_{K}^{2}}, n_{M}^{2}n_{K}^{2}\right)$. Since $N \leq_{(m,\mathcal{ES})} K$,
then, there will be two cases:\\

i. If $m_{N} < m_{K}$, then, \\
\begin{eqnarray*}
m_{M +^{P} N}^{2} - m_{M +^{P} K}^{2} = (m_{M}^{2}- m_{K}^{2})(1 - m_{M}^{2}) <0
\end{eqnarray*}

then, $m_{M +^{P} N} < m_{M +^{P} K}$.\\

ii. Consider the $m_{N} = m_{K}$ and $\mathcal{ES}_{N} \leq \mathcal{ES}_{K}$. Then, we have
\begin{eqnarray*}
m_{M +^{P} N}^{2} =m_{M}^{2}+m_{N}^{2}-m_{M}^{2}m_{N}^{2}=m_{M}^{2}+m_{K}^{2}-m_{M}^{2}m_{K}^{2} = m_{M +^{P}K}^{2},
\end{eqnarray*}
and
\begin{eqnarray*}
n_{N}^{2} = m_{N}^{2} +1 - 2\mathcal{ES}_{N} \geq m_{K}^{2} +1 - 2\mathcal{ES}_{K} = n_{K}^{2}.
\end{eqnarray*}
Then,
\begin{eqnarray*}
n_{M +^{P} N}^{2} = n_{M}^{2} n_{N}^{2}\geq n_{M}^{2}n_{K}^{2}= n_{M +^{P} K}^{2}.
\end{eqnarray*}
It follows from that
\begin{eqnarray*}
\mathcal{ES}_{M +^{P} N}=\frac{m_{M +^{P} N}^{2}- n_{M +^{P} N}^{2}+1}{2} \geq \frac{m_{M +^{P} K}^{2}- n_{M +^{P} K}^{2}+1}{2}=\mathcal{ES}_{M +^{P} K}.
\end{eqnarray*}
From the cases (i) and (ii), we can write $M +^{P} N \leq_{(m,\mathcal{ES})} M +s^{P} K$.
\end{proof}

\begin{Corollary}
Let For $M=(m_{M}, n_{M}), N=(m_{N}, n_{N})\in \mathcal{L}$ and $\alpha, \alpha_{1}, \alpha_{2}\in \mathbb{R}^{+}$.
Then,
\begin{itemize}
\item [i.] If $M \leq_{(m,\mathcal{ES})} N $, then $\alpha M \leq_{(m,\mathcal{ES})} \alpha N $,
\item[ii.] If $\alpha_{1} \leq \alpha_{2}$, then $\alpha_{1} M \leq_{(m,\mathcal{ES})} \alpha_{2} M $.
\end{itemize}
\end{Corollary}

\subsection{Algorithm}

Pythagorean fuzzy weighted averaging operator(PFWA) was given by Yager \cite{Yager1}.
The following definition is based on the definition of PFWA of Yager.

\begin{Definition}
Let $\pi_{i}=(m_{i}, n_{i})$ be PFVs in $\mathcal{L}$, $\omega=(\omega_{1}, \omega_{2},\cdots, \omega_{k})^{T}$
is the weighted vector such that for $i=1,2,\cdots,k$, $\omega_{i}\in [0,1]$ with
$\sum_{i=1}^{k}\omega_{i}=1$. Then the mapping $\Phi_{PFWA}=\mathcal{L}^{n}\rightarrow \mathcal{L}$ given by
\begin{eqnarray*}
\Phi_{PFWA}(\pi_{i})=\left(\sum_{i=1}^{k}\omega_{i}m_{i}, \sum_{i=1}^{k}\omega_{i}n_{i}\right)
\end{eqnarray*}
is called PFWA opeartor.
\end{Definition}

\begin{Definition}\cite{garg}, \cite{Yager1}, \cite{YagerAbb} \label{def:90}
The mapping $\mathcal{PF}: \mathcal{L}^{k} \rightarrow \mathcal{L}$ given by
\begin{eqnarray*}
\mathcal{PF}_{\omega}(\pi_{i})=\omega_{1}\pi_{1} +_{P} \omega_{2}\pi_{2} +_{P} \cdots +_{P} \omega_{k}\pi_{k} = \left( \sqrt{1-\prod_{i=1}^{k}(1-m_{i}^{2})^{\omega_{i}}}, \prod_{i=1}^{k}n_{i}^{\omega_{i}}\right)
\end{eqnarray*}
is called the $k$ dimensional PFWA operator.
\end{Definition}

Definition \ref{def:90} can be used to simplify the computation concerning PFWA operators.

\begin{Definition}
Let $\varphi_{X} \in \Phi(U)$. Then,
\begin{eqnarray*}
\mathcal{AP}_{\varphi_{X}}(a)=\biguplus_{x\in X}\frac{\mathcal{ES}_{f_{X}(x)}}{\sum_{x\in X}\mathcal{ES}_{f_{X}}(x)}\varphi(x)(a)
\end{eqnarray*}
is called the aggregated Pythagorean fuzzy decision value(APFDV) of $a\in \varphi_{X}$.
\end{Definition}

\textbf{Algorithm:}\\

\begin{itemize}
\item \textbf{Step 1:} Let $U=\{p_{1}, p_{2}, \cdots, p_{i}\}$ and $Z=X\cap Y=\{s_{1}, s_{2}, \cdots, s_{j}\}$.
Choose PFSSs $\varphi_{X}, \varphi_{Y}$ over $U$ which are separately two groups with $f_{X}, f_{Y}$.
\item \textbf{Step 2:} Calculate the extended intersection $\varphi_{X} \hat{\cap}_{E} \varphi_{Y}=\varphi_{Z}$.
\item \textbf{Step 3:} For $k=1,2,\cdots,i$, calculate the APFDVs
\begin{eqnarray*}
\mathcal{AP}_{\varphi_{Z}(p_{k})}=\biguplus_{\ell=1}^{j}\frac{\mathcal{ES}_{f_{Z}(x_{\ell})}}{\sum_{\ell=1}^{j}\mathcal{ES}_{f_{Z}}(x_{\ell})}\varphi(Z)(x)(p_{k})
\end{eqnarray*}
\item \textbf{Step 4:} Rank $\mathcal{AP}_{\varphi_{Z}(p_{k})}$, $(k=1,2,\cdots,i)$ descendingly under the order $\leq_{(m, \mathcal{ES})}$.
\item \textbf{Step 5:} Rank $p_{j}$, $(j=1,2,\cdots,k)$ correspondingly and output $p_{k}$ as the optimal decision, if
$\mathcal{AP}_{\varphi_{Z}(p_{i})}$ is the largest PFV under the order $\leq_{(m, \mathcal{ES})}$.
\end{itemize}

\subsection{Application}

We consider the values of Table \ref{table:1}, Table \ref{table:2} and Table \ref{table:4}.
From the $f_{Z}$, we compute the "expectation values"
$\mathcal{ES}_{f_{Z}(x_{\ell})}$, that reveal the weight vector
\begin{eqnarray*}
\omega=\{0.21001927, 0.12524085, 0.27938343, 0.14065510, 0.24470135\}^{T}
\end{eqnarray*}
to be used for calculating the APFDVs. The weight vector $\omega$ computed as $\frac{\mathcal{ES}_{f_{Z}(x_{\ell})}}{\sum_{\ell=1}^{j}\mathcal{ES}_{f_{Z}}(x_{\ell})}$ (\ref{table:7}). The $\mathcal{AP}_{\varphi_{Z}(p_{k})}$ in Step 3 is found as
\begin{eqnarray*}
\mathcal{AP}_{\varphi_{Z}(p_{i})}=\mathcal{PF}_{\omega}(\pi_{i})\bigg( \varphi_{Z}(s_{1})(p_{i}), \varphi_{Z}(s_{2})(p_{i}), \varphi_{Z}(s_{3})(p_{i}), \varphi_{Z}(s_{5})(p_{i}), \varphi_{Z}(s_{6})(p_{i})\bigg).
\end{eqnarray*}
For example, $\mathcal{AP}_{\varphi_{Z}(p_{1})}=(0.6314, 0.6434)$.\\

From Table \ref{table:8}, we have,
\begin{eqnarray*}
\mathcal{AP}_{\varphi_{Z}(p_{2})} \leq_{(m, \mathcal{ES})} \mathcal{AP}_{\varphi_{Z}(p_{3})} \leq_{(m, \mathcal{ES})} \mathcal{AP}_{\varphi_{Z}(p_{4})} \leq_{(m, \mathcal{ES})} \mathcal{AP}_{\varphi_{Z}(p_{1})}.
\end{eqnarray*}

According to these results, the patients will be sorted as follows:

\begin{eqnarray*}
p_{4} > p_{3} > p_{1} > p_{2}.
\end{eqnarray*}

\begin{table}
 \caption{}\label{table:7}
  \vglue2mm
\centering
 {
\begin{tabular}{c  c  c  c  c c}
\hline
S   & $s_{1}$ & $s_{2}$ & $s_{3}$ & $s_{5}$ & $s_{6}$ \\
\hline
$f_{Z}$  & (0.5, 0.4) & (0.1, 0.6) & (0.7, 0.2) & (0.3, 0.6) & (0.6, 0.3)\\
$\mathcal{ES}_{f_{Z}(x_{\ell})}$ & 0.545 & 0.325 & 0.725 & 0.365 & 0.635 \\
$\omega_{\ell}$  & 0.2100193 & 0.1252409 & 0.2793834 & 0.1406551 & 0.2447013\\
\hline
\end{tabular}}
\end{table}

\begin{table}
 \caption{Measures}\label{table:8}
  \vglue2mm
\centering
 {
\begin{tabular}{c  c  c  c  c }
\hline
P   & $APFDVs$ & $\mathcal{ES}(\mathcal{AP}_{\varphi_{Z}(p_{k})})$ & $\mathcal{SF}(\mathcal{AP}_{\varphi_{Z}(p_{k})})$ & $\mathcal{AF}(\mathcal{AP}_{\varphi_{Z}(p_{k})})$ \\
\hline
$p_{1}$  & (0.6314, 0.6434) & 0.4923512 & -0.0152976 & 0.81262952 \\
$p_{2}$ & (0.,3601, 0.5271) & 0.4259188 & -0.1481624 & 0.40750642  \\
$p_{3}$  & (0.5156, 0.4358) & 0.53796086 & 0.07592172 & 0.455765 \\
$p_{4}$  & (0.5554, 0.3642) & 0.58791376 & 0.17582752 & 0.4411108 \\
\hline
\end{tabular}}
\end{table}

\section{Conclusion}
Pythagorean fuzzy sets was initiated by Yager \cite{Yager0}.
Many researchers are concerned with decision-making issues in the PFS environment \cite{garg, Guleria,Pengetal, PengYang, PengYang2, YagerAbb, Yager1, YagerKit}.
PFSSs were given by Peng et al. \cite{Pengetal}. In this paper, we introduce the definition of new PFSS with a parameter.
The properties and some operations of the new SS are examined.
For the decision-making process, the functions of score, accuracy and expectation score are defined.
Depending on PFWA, the formula of aggregated Pythagorean fuzzy decision value was obtained and given an algorithm.

\end{document}